\documentclass[preprint]{elsarticle}

\usepackage{amssymb,amsmath,amsthm}

\usepackage[dvips]{epsfig}
\usepackage{epsf}

\advance \topmargin by -\headheight
\advance \topmargin by -\headsep
\evensidemargin \oddsidemargin
\def \proof{\bigbreak\noindent{\it Proof.\ \ }}

\def \endpf{{\ \ $\Box$ \medbreak}}

\newtheorem{lemma}{Lemma}[section]

\newtheorem{proposition}[lemma]{Proposition}
\newtheorem{theorem}[lemma]{Theorem}

\theoremstyle{definition}
\newtheorem{definition}[lemma]{Definition}

\newtheorem{remark}[lemma]{Remark}

\begin{document}

\begin{frontmatter}
\title{Inconstancy of finite and infinite sequences}

\author[jus]{Jean-Paul Allouche\corref{cor1}}
\ead{allouche@math.jussieu.fr}

\author[rte]{Laurence Maillard-Teyssier\corref{cor2}}
\ead{Laurence.Teyssier-Maillard@rte-france.com}

\address[jus]{CNRS, Institut de Math.,
Universit\'e P. et M. Curie, Case 189,
4 Place Jussieu, F-75252 Paris Cedex 05, France}

\address[rte]{RTE, DMA, Immeuble Le Colbert, 9 rue de la Porte de Buc, 
BP 561, 78005 Versailles Cedex, France} 

\cortext[cor1]{Corresponding author}

\begin{abstract}
In order to study large variations or fluctuations of finite or infinite
sequences (time series), we bring to light an 1868 paper of Crofton and
the (Cauchy-)Crofton theorem. After surveying occurrences of this result
in the literature, we introduce the {\em inconstancy\,} of a sequence and 
we show why it seems more pertinent than other criteria for measuring its
variational complexity. We also compute the inconstancy of classical binary 
sequences including some automatic sequences and Sturmian sequences.
\end{abstract}

\begin{keyword}Fluctuations \sep time series \sep discrete curves \sep
Cauchy-Crofton theorem \sep inconstancy of sequences \sep entropy 
\sep automatic sequences \sep Sturmian sequences \sep combinatorics on words

\MSC[2000] 53C65 \sep 52C45 \sep 37M10 \sep 11B85 \sep 52A38 \sep
           51M25 \sep 28A75 \sep 62M10 \sep 52A10 \sep 68R15
\end{keyword}

\end{frontmatter}

\medskip

\begin{flushright}
{\it The voyage of the best ship is a zig-zag line of a hundred tacks. \\
See the line from a sufficient distance, and it straightens itself to \\
the average tendency... (Ralph Waldo Emerson, Emerson Essays, 1899)}
\end{flushright}

\section{Introduction}

How is it possible to define {\em and\,} to detect large variations or 
fluctuations of a sequence (with possible applications to the [discrete] 
time evolution of biological, financial, musical phenomena and so on). The 
usual approach is based on computing the distance of the associated 
{\em piecewise affine function} to the corresponding linear regression 
line, i.e., on computing the {\em residual variance}. But this quantity 
somehow describes total distance to ``regularity'', and says nothing 
about possibly large local fluctuations: for example, it may not 
discriminate between an exponentially growing function and a fractal-like 
``chaotic'' (disordered) curve. In particular one should remember that 
dictionaries defining ``fluctuation'' use words with a similar meaning 
among which ``wavering'', ``unsteadiness'', ``vacillation'', ``erraticness'', 
``variability'', etc.

\medskip

We suggest here to bring to light -- especially for applications to sequences
-- a paper of Crofton dated 1868 \cite{Crofton} (see also the papers of 
Cauchy \cite{Cauchy1, Cauchy2} and the papers of Steinhaus \cite{Steinhaus} and 
of Dupain, Kamae and Mend\`es France \cite{DKMF}). Crofton studies the average 
number of intersection points of a curve with random straight lines. But this average 
number can be thought of as a measure of the fluctuations of the curve.
Namely, for a straight line or a curve ``looking like a straight line'', this average 
number is equal to $1$, while it has a very large value for a ``very complicated''
curve. Following this idea, we propose a measure of large variations of a sequence
and we compare it with the residual variance. Conversely, this measure will 
allow us to {\em decide\,} whether a sequence is ``more complicated'' than another 
in cases where the visual aspect does not suffice to suggest an intuitive answer. 
We will also show that this measure can be applied to {\em infinite sequences\,} 
satisfying some technical condition (in particular certain {\em automatic 
sequences\,} as well as {\em Sturmian sequences}; see, e.g., \cite{AS2}) to 
describe their ``complexity''.

\medskip

As will be recalled, the ideas of Cauchy and Crofton were already used in various 
contexts: one of our purposes is {\em to insist on their usefulness for measuring 
the complexity of discrete phenomena}, as a compromise between measuring intensity, 
time and consecutive repetitions. These ideas will be applied in a subsequent paper 
(see \cite{VMAH}) to fluctuations of biological parameters, e.g., the weight, or the 
{\it Quetelet index}\footnote{In \cite{Quetelet} Quetelet asserts that weights
vary like heights squared for adults but more like (heights)$^{5/2}$ for children
(see p. 52--53, and p. 61), while the ``simplified'' definition of the BMI is the ratio
of the weight by the height squared.}, often called the {\it BMI\,} (Body Mass Index; see,
e.g., \cite{Quetelet, BKT, RC82, RC84, RC2005}) for children: are ``large fluctuations'' 
of the BMI risk factors for cardiovascular diseases in relation with 
the {\em metabolic syndrome}? This question was addressed with other tools in 
\cite{Vergnaud} (see also the references therein). 
We also aim to try to apply this measure of fluctuations to other questions, e.g., 
analyzing fluctuations of the stockmarket, and quantifying the ``smoothness'' of 
musical themes. 

\section{Defining the {\em Inconstancy\,} of a curve}

A possible approach for describing large variations or large fluctuations of a 
curve is to ``compare'' it with a straight line. More precisely we can count the 
number of intersection points of random straight lines with the given curve: if 
this number is small on average, the curve behaves roughly as a straight line; if 
this number is large, the curve is ``complicated''. Is there an ``easy'' way 
to compute this number? The Cauchy-Crofton theorem answers the question.

\subsection{The Cauchy-Crofton theorem}

Consider a plane curve $\Gamma$. Let $\ell(\Gamma)$ denote the length of $\Gamma$ 
and let $\delta(\Gamma)$ denote the perimeter of the closed curve forming the 
edge of the convex hull of $\Gamma$.
Let $\Omega(\Gamma)$ be the set of straight lines which intersect $\Gamma$. 
Any line can be defined as the set of $(x,y)$ such that 
$x \cos\theta + y \sin\theta - \rho = 0$, where $\theta$ belongs to $[0,\pi)$ 
and $\rho$ is a positive real number. A straight line is therefore completely determined
by $(\rho, \theta)$. Letting $\mu$ denote the Lebesgue measure on the set
$\{(\rho, \theta), \ \rho \geq 0, \ \theta \in [0,\pi)\}$, the {\it average\,} 
number of intersection points between the curve $\Gamma$ and a line in $\Omega$ is
defined by
$$
{\mathcal N}(\Gamma) := \int_{D \in \Omega(\Gamma)} \sharp(\Gamma \cap D) 
\frac{{\rm d}\rho \ {\rm d}\theta}{\mu(\Omega(\Gamma))}\cdot
$$
The following result can be found in \cite[p.~184--185]{Crofton}; see also the
papers of Cauchy \cite{Cauchy1, Cauchy2}.

\begin{theorem}[Cauchy-Crofton]\label{Crofton}
The average number of intersection points between the curve $\Gamma$ and the 
straight lines in $\Omega$ satisfies
$$
{\mathcal N}(\Gamma) = \frac{2 \ell(\Gamma)}{\delta(\Gamma)}\cdot
$$
\end{theorem}

\begin{remark}
In his paper, Crofton speaks of ``Local or Geometrical Probability''; he writes
about Probabilities, ``The rigorous precision, as well as the extreme beauty of the
methods and results... the subtlety and delicacy of the reasoning...'', and he
quotes Laplace: ``ce calcul d\'elicat''. 
Crofton's result is explained in Steinhaus' paper \cite{Steinhaus}. It is presented 
in an illuminating way with several examples in the paper of Dupain, Kamae, and 
Mend\`es France \cite{DKMF}: these authors studied the notion of {\em entropy 
of a curve} and of {\it temperature of a curve} introduced by Mend\`es France in
\cite{MF83}. Note that the occurrence of the number $2$ in the numerator can be 
understood by considering the case where $\Gamma$ is a segment: the average number 
of intersection points is equal to $1$, while the perimeter of the convex hull of the 
segment is twice the length of the segment (why twice? go back to the definition
[``closed curve...''] or think of the case where the segment is replaced by a thin 
rectangle whose width tends to zero).
\end{remark}

\begin{remark}
The reader will have noted that Crofton's approach has much to do with
the famous Buffon needle problem \cite[p.~100--104]{Buffon}
also known as the Buffon-Laplace needle problem; see \cite[p.~359--360]{Laplace}.
The area of this type of result is known as ``Integral Geometry''. This 
terminology seems to have been introduced by Blaschke in his ``Vorlesungen \"uber 
Integralgeometrie'' \cite{Blaschke1, Blaschke2}. More recent references are 
the book of Santal\'o \cite{Santalo}, and the forthcoming book of Langevin 
\cite{Langevin2} (see also \cite{Langevin1}). An interesting review of the books 
of Blaschke and of the 1936 edition of the book of Santal\'o is \cite{Myers}. A 
nice exposition of the (proof of the) theorem of Cauchy-Crofton, where the curve is 
only supposed to be rectifiable, can be found in the paper of Ayari and Dubuc \cite{AD}. 
We also recommend for a first approach the texts of Mend\`es France \cite{MF06} 
and of Teissier \cite{Teissier}. 
Note that the Crofton theorem is also (and more correctly) called the Cauchy-Crofton
theorem in the literature.
\end{remark}

\begin{remark}
Using the theorem of Cauchy-Crofton to define a measure of complexity of a curve was 
first suggested by Mend\`es France \cite[page~92]{MF82}; also see \cite{MF83a, MF83}. 
It was also proposed later, e.g., in \cite{CYW} where the name ``folding index'' is used. 
Also note that the Crofton formulas in \cite{Crofton} are used frequently in many fields. These include complex motor behaviour in human movements \cite{Cor1} (also 
see \cite{Cor2, Cor3}), study of human blood and transfusion \cite{WVDGHKLMN}, simulation 
of gravitational evolution \cite{SBMSSS}, anisotropies of the secondary cosmic microwave 
background \cite{GS}, grain size distribution analysis for polycrystalline thin films 
\cite{CMSPNZCA}, image analysis of crystalline agglomerates \cite{PPVA}, measurement of 
convolution in cotton fibers \cite{HCLB}, all applications of LIS (Line-Intercept Sampling), 
e.g., to the statistical analysis of vegetation or wildlife, see for example \cite{ZPN} and 
the references therein (in particular \cite{Kaiser} in the references below), spatial analysis 
of urban maps \cite{Furuyama}, in a discussion about examples of information processing 
coming from neurophysiology, cognitive psychology, and perception 
\cite[pp.~1182--1185]{Nicolis}, and even relations between art and complexity 
\cite{MFH} (also see \cite{Nesetril1, Nesetril2}). 
\end{remark}

\subsection{The {\em inconstancy} of a curve}

The theorem of Cauchy-Crofton suggests the following definition.

\begin{definition}
Let $\Gamma$ be a plane curve of length $\ell(\Gamma)$ and such that 
the perimeter of its convex hull is equal to $\delta(\Gamma)$.
The {\em inconstancy\,} of the curve $\Gamma$, denoted ${\mathcal I}(\Gamma)$, 
is defined by
$$
{\mathcal I}(\Gamma) := \frac{2 \ell(\Gamma)}{\delta(\Gamma)}\cdot 
$$

\end{definition}

\begin{remark}
The above definition and the Cauchy-Crofton theorem show that the inconstancy
of the union of two curves is at most the sum of the inconstancies of these curves,
that the inconstancy of a curve is equal to the inconstancy of its translated,
rotated or homothetic curve, etc.
\end{remark}

\section{Comparison with other criteria}

Other criteria for measuring fluctuations of a discrete curve can be found in the
literature for real (e.g., biological) phenomena: qualitative classification with 
predetermined cut-off points, maximal values, residual variance, etc. (see, e.g., 
the discussion in \cite[pp.\ 316--317]{Vergnaud} for weight fluctuations). By 
oversimplifying most of the various definitions, one could say that they aim to measure
the ``distance'' between the considered curve and a straight line, but this distance
can be computed globally or locally. We recall the definition of {\em regression
line}, of {\em residual variance}, and of {\em mean square error}.

\begin{definition}\label{def}
Let $(x_i, y_i)_{i = 1,2, \ldots, n}$ be a family of $n \geq 3$ points. 
Their {\em regression line} is the straight line that minimizes the sum of squares 
of distances from the $(x_i, y_i)$'s to it. 
Letting $\overline{x} = (\sum_{1 \leq i \leq n} x_i)/n$ and
$\overline{y} = (\sum_{1 \leq i \leq n} y_i)/n$ denote the averages of 
the $x_i$'s and of the $y_i$'s, the equation of the regression line is
$$
Y = \hat{a} X + \hat{b}, \ \mbox{\rm where }
\hat{a} = \frac{\sum_{1 \leq i \leq n}(x_i - \overline{x})(y_i - \overline{y})}
               {\sum_{1 \leq i \leq n}(x_i-\overline{x})^2} \ \ 
               \mbox{\rm is the {\it correlation}, and } \
\hat{b} = \overline{y} - \hat{a} \, \overline{x}.
$$
The {\em MSE} (i.e., {\em mean square error}) and the {\em RMSE} 
(i.e., {\em root mean square error}) of the $(x_i, y_i)$'s are defined by
$$
MSE : = \frac{1}{n-2} \sum_{1 \leq i \leq n} (y_i - \hat{a} \, x_i - \hat{b})^2 \ \
\mbox{and } \ RMSE := \sqrt{MSE}. 
$$
The mean square error is sometimes called {\em residual variance}.

\end{definition}

\bigskip

We also introduce some notation.

\begin{definition}
Let $n$ be a positive integer. We define $\Gamma(a_1, a_2, \ldots, a_n)$ to 
be the union of the $n$ segments $(0,0)$---$(1, a_1)$, $(1, a_1)$---$(2, a_2)$, 
$\ldots$ $(n-1, a_{n-1})$---$(n, a_n)$. (Note that we have $(n+1)$ points, and
that, without loss of generality, we suppose that the curve begins at the origin.)
\end{definition}

\subsection{Why is MSE not satisfactory to measure fluctuations?}

\begin{figure}[h]
\begin{center}
        \epsfig{file=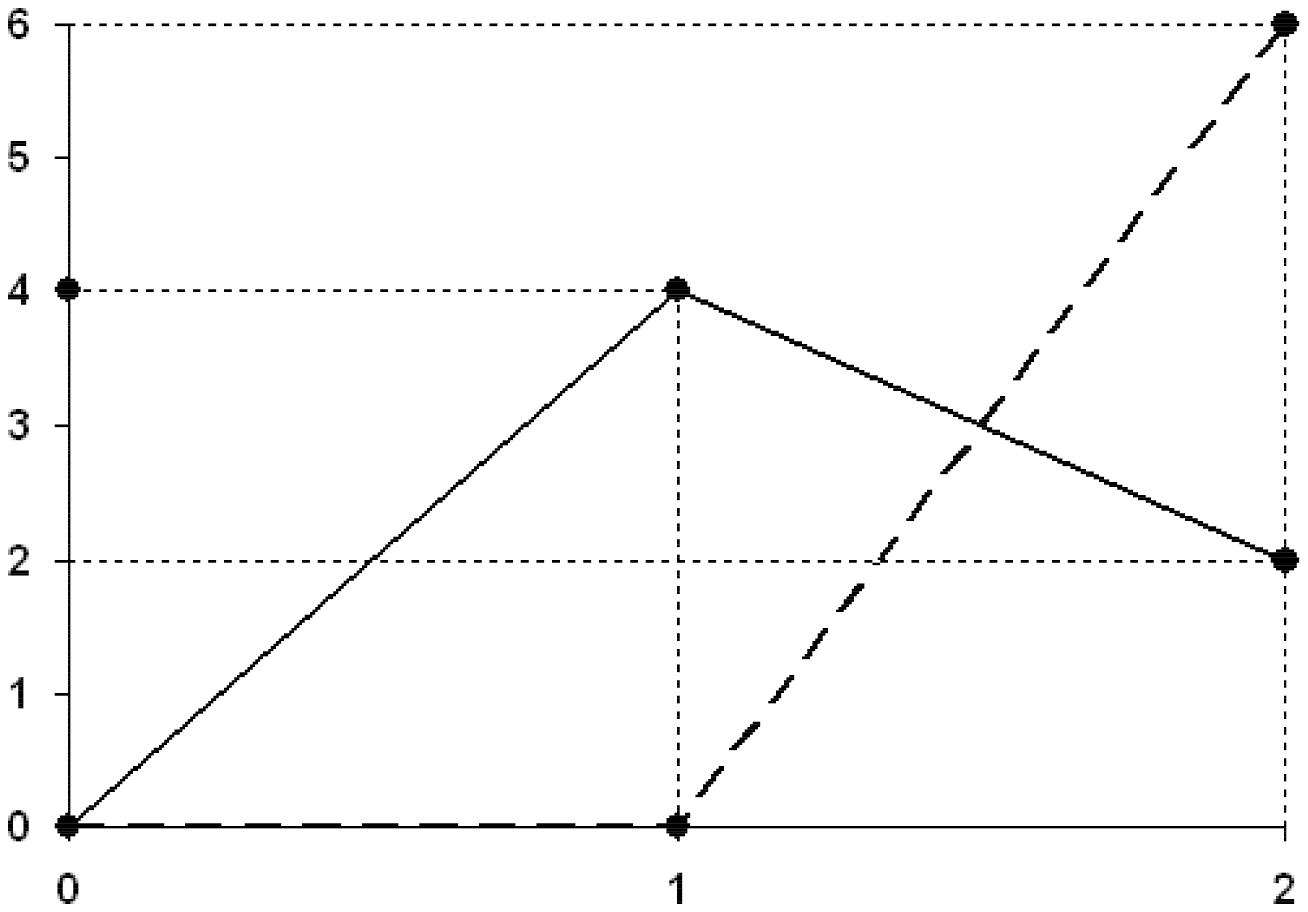}
\end{center}
\caption{Same MSE}
\label{same-mse}
\end{figure}

In this section we show two curves having same length: one is ``fluctuating'', 
the other increases quickly, but their residual variances are both equal to $6$;
see Figure~\ref{same-mse}. Note that when we say that the first curve is more
``fluctuating'' than the second one, it means for example that for a variation
of weight or of BMI, the first curve is really fluctuating, while the second
one just shows some (possibly quick) growth (also see Remarks~\ref{growth}
and the beginning of Section~\ref{nonintuitive} below).

\subsection{Comparing MSE and inconstancy}

Are residual variance and inconstancy of a curve comparable? We will
prove that this is not the case, even for very simple curves, thanks to
two easy lemmas.

\begin{lemma}\label{res}
Let ${\mathcal R}(\Gamma(a_1, a_2))$ be the residual variance of the
curve $\Gamma(a_1, a_2)$. Then 
$$
{\mathcal R}(\Gamma(a_1, a_2)) = \frac{(2a_1 - a_2)^2}{6} \cdot
$$
\end{lemma}

\proof
The computation is straightforward. The linear regression straight line
is parallel to $(0,0)$---$(2, a_2)$, and it contains the center of gravity
of the triangle  $(0,0)$, $(1, a_1)$, $(2, a_2)$. Or simply compute from
Definition~\ref{def}: $\overline{x} = 1$, $\overline{y} = (a_1 + a_2)/3$,
$\hat{a} = a_2/2$, and $\hat{b} = (2a_1 - a_2)/6$, hence 
${\mathcal R}(\Gamma(a_1, a_2)) = (2a_1 - a_2)^2/6$. \endpf

\begin{lemma}\label{inc}
Let $\Gamma(a_1, a_2)$ be the curve defined as the union of the two straight
line segments $(0,0)$---$(1, a_1)$ and $(1, a_1)$---$(2, a_2)$. Then,
${\mathcal I}(\Gamma(a_1, a_2))$, the inconstancy of $\Gamma(a_1, a_2)$,
is given by
$$
{\mathcal I}(\Gamma(a_1, a_2)) = 
\dfrac{2}{1+\dfrac{\sqrt{{a_2}^2 + 4}}{\sqrt{{a_1}^2+1} + \sqrt{(a_2-a_1)^2 + 1}}}
$$
\end{lemma}

\proof
The proof is again straightforward. The length of $\Gamma(a_1, a_2)$ and the
perimeter of the convex hull of $\Gamma(a_1, a_2)$ are given respectively by
$$
\sqrt{a_1^2+1} + \sqrt{(a_2-a_1)^2+1} 
\ \ \mbox{\rm and } \
\sqrt{a_1^2+1} + \sqrt{(a_2-a_1)^2+1} + \sqrt{a_2^2 + 4}. 
\ \ \ \ \ \ \ \ \ \ \ \ \ \ \Box 
$$

We can now state the non-comparability of residual variance and inconstancy.

\begin{figure}[h]
\begin{center}
        \epsfig{file=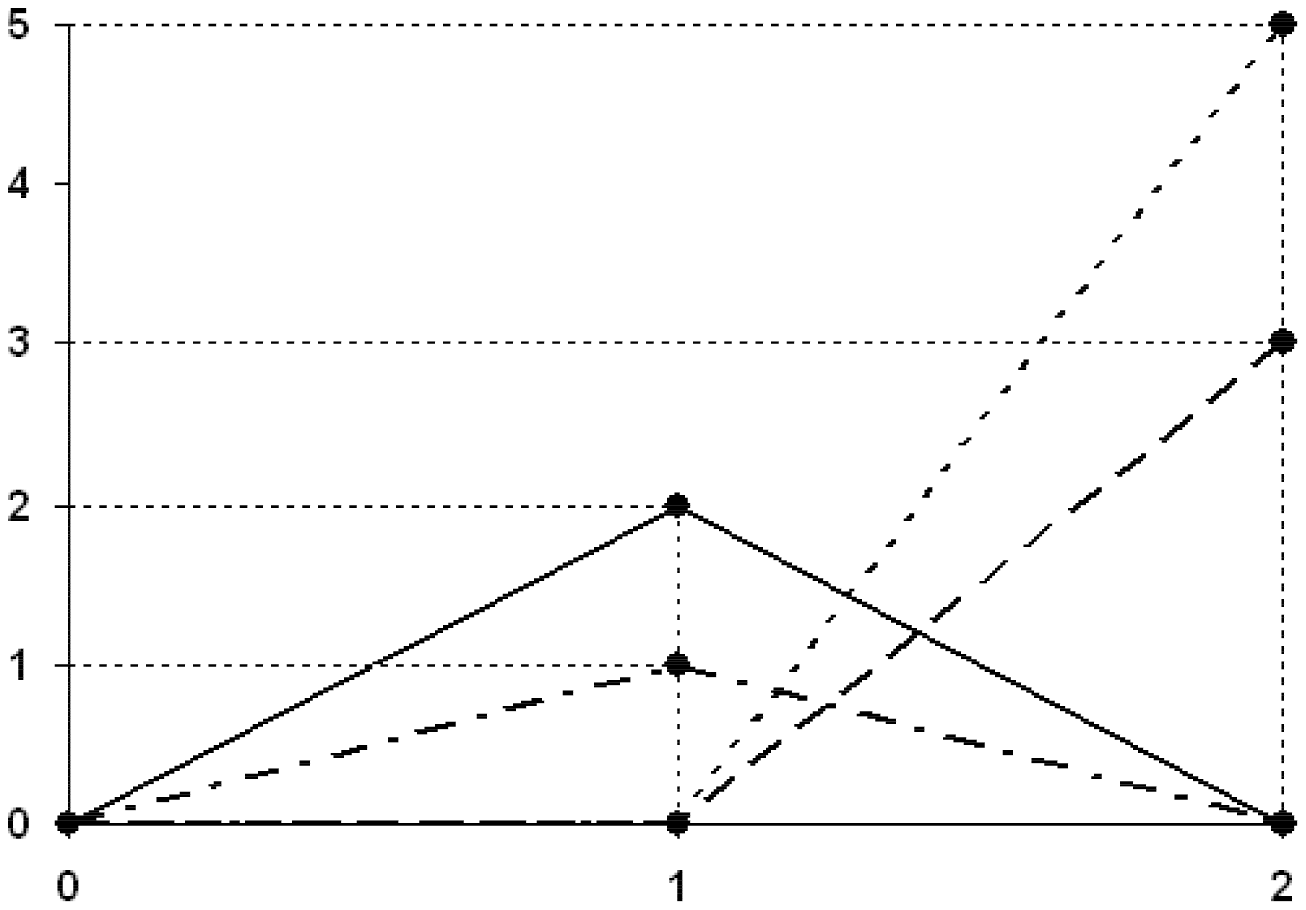}
\end{center}
\caption{Comparing residual variance and inconstancy}
\label{proposition}
\end{figure}

\begin{proposition}\label{compare}
Residual variance and inconstancy of a curve are not comparable.
More precisely, there exist four curves $\Gamma_i$, $i=1, 2, 3, 4$,
(see Figure~\ref{proposition} and the proof below) such that, if 
${\mathcal R}(\Gamma_i)$ and ${\mathcal I}(\Gamma_i)$ are their 
residual variances and inconstancies, then the following inequalities hold:
$$
\begin{array}{llllllll}
&{\mathcal R}(\Gamma_1) &<&{\mathcal R}(\Gamma_2)  
&<&{\mathcal R}(\Gamma_3) &< &{\mathcal R}(\Gamma_4) \\
&{\mathcal I}(\Gamma_4) &<&{\mathcal I}(\Gamma_2) 
&< &{\mathcal I}(\Gamma_1) &< &{\mathcal I}(\Gamma_3).
\end{array}
$$
\end{proposition}

\proof Using Lemmas~\ref{res} and \ref{inc} above, we get the residual 
variances ${\mathcal R}(\Gamma_i)$ and inconstancies ${\mathcal I}(\Gamma_i)$ 
of the following curves $\Gamma_i$
$$
\begin{array}{llllllllll}
\Gamma_1 &:=& \Gamma(1,0) \qquad {\mathcal R}(\Gamma_1) &=& \frac{2}{3}
            &\approx 0.67 \qquad {\mathcal I}(\Gamma_1) &=& \frac{2\sqrt{2}}{1+\sqrt{2}}
                                &\approx& 1.17 \\
\Gamma_2 &:=& \Gamma(0,3) \qquad {\mathcal R}(\Gamma_2) &=& \frac{3}{2}
            &\approx 1.50 \qquad {\mathcal I}(\Gamma_2) &=& \frac{2+2\sqrt{10}}{1+\sqrt{10}+\sqrt{13}}
                                &\approx& 1.07  \\
\Gamma_3 &:=& \Gamma(2,0) \qquad {\mathcal R}(\Gamma_3) &=& \frac{8}{3}
            &\approx 2.67 \qquad {\mathcal I}(\Gamma_3) &=& \frac{2\sqrt{5}}{1 + \sqrt 5}   
                                &\approx& 1.38  \\
\Gamma_4 &:=& \Gamma(0,5) \qquad {\mathcal R}(\Gamma_4) &=& \frac{25}{6}
            &\approx 4.17 \qquad {\mathcal I}(\Gamma_4) &=& \frac{2+2\sqrt{26}}{1+\sqrt{26}+\sqrt{29}}
                                &\approx& 1.06  
\ \ \ \ \ \ \ \ \ \ \ \ \ \ \ \Box \\
\end{array} 
$$

\begin{remark}\label{growth}
Comparing ${\mathcal I}(\Gamma_2)$ and ${\mathcal I}(\Gamma_4)$ shows again that
``fluctuating'' is not the same as ``growing''. More generally, with the notation
of Lemma~\ref{inc} above, looking at ${\mathcal I}(\Gamma(0,x))$, shows
that ${\mathcal I}(\Gamma(0,0)) = 1 = \lim_{x \to \infty}{\mathcal I}(\Gamma(0,x))$.
When $x$ varies from $0$ to $\infty$ the quantity ${\mathcal I}(\Gamma(0,x))$
increases from $1$ to a small value $>1$ then it decreases back to $1$.
\end{remark}

\begin{remark}
There are other quantities that also ``measure'' the fluctuations of a curve.
For example, keeping the notations of Definition~\ref{def}: the {\em total variation}
is defined as the mean of $(y_i - \overline{y})^2$, i.e., as 
$(\sum_{1 \leq i \leq n} (y_i - \overline{y})^2)/n$; the {\em maximal distance} is 
defined as $\max_{1 \leq i \leq n} |y_i - \hat{a} \, x_i - \hat{b}|$.
The reader can easily compute these quantities for the curve $\Gamma(a_1, a_2)$ and 
check that they are not comparable to the inconstancy of $\Gamma(a_1, a_2)$.  

\medskip

Total variation: $\dfrac{2(a_1^2 + a_2^2 - a_1a_2)}{9}\cdot$ \ \ \ \ \ 
Maximal distance: $\dfrac{|2a_1 - a_2|}{3}\cdot$ 
\end{remark}

\section{Pertinence of the use of inconstancy: simple arguments}

\subsection{A single fluctuation}

\begin{figure}[h]
\begin{center}
         \epsfig{file=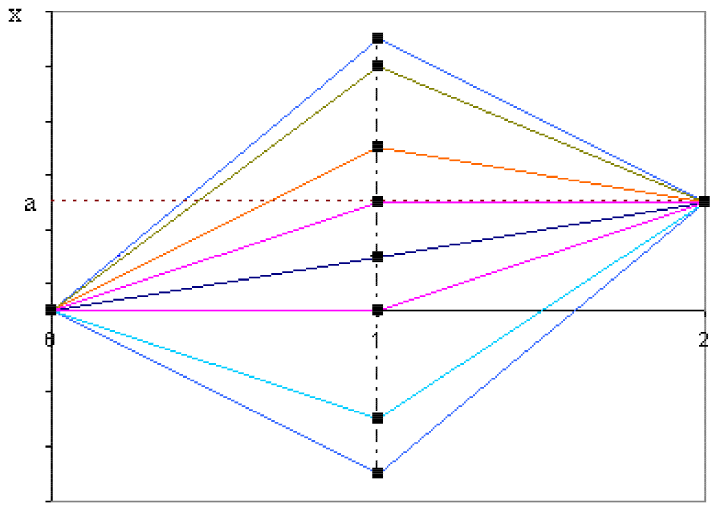}
\end{center}
\caption{Varying the intermediate value}
\label{3varyingpoints}
\end{figure}

Taking again the example in the previous section of a curve consisting of
two straight line segments, let us vary the value $a_1$, say $x := a_1$, and
fix $a_2 = a$ (see Figure~\ref{3varyingpoints}). The inconstancy 
${\mathcal I}(\Gamma(x, a))$ is thus given by
$$
{\mathcal I}(\Gamma(x, a)) =
\dfrac{2}{1+\dfrac{\sqrt{{a}^2 + 4}}{\sqrt{{x}^2+1} + \sqrt{(a-x)^2 + 1}}}
$$
This map $x \to {\mathcal I}(\Gamma(x, a))$ is increasing for $x \geq a/2$, 
which is in agreement with what a ``fluctuation'' should be.

\begin{figure}[h]
\begin{center}
        \epsfig{file=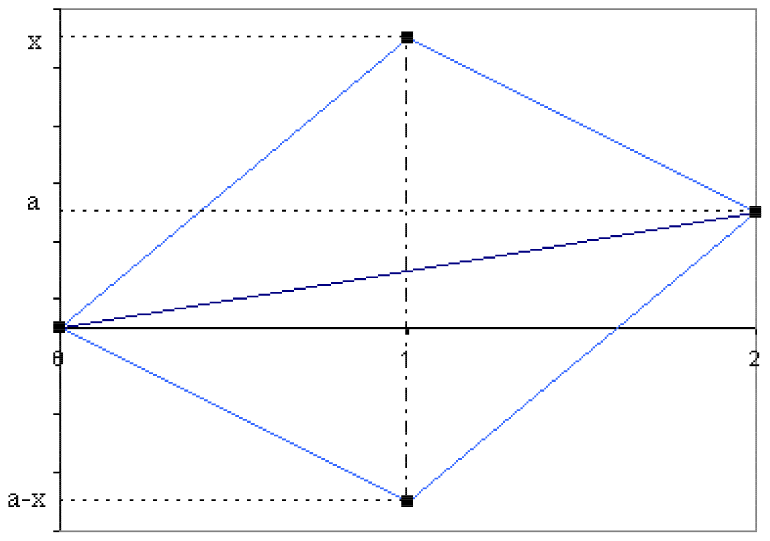}
\end{center}
\caption{Symmetry}
\label{symmetry}
\end{figure}

It is clear that ${\mathcal I}(\Gamma(x, a)) = {\mathcal I}(\Gamma(a-x, a))$, which
shows that the line $x = a/2$ is a symmetry axis. In other words, ``exchanging''
the two segments, more precisely replacing ($(0,0)$---$(1,x)$), ($(1,x)$---$(2,a)$)
by ($(0,0)$---$(1,a-x)$), ($(1,a-x)$---$(2,a)$), does not change the inconstancy
(see Figure~\ref{symmetry}). Of course this is a necessary condition for a fluctuation
criterion.

\begin{figure}[!h]
\begin{center}
        \epsfig{file=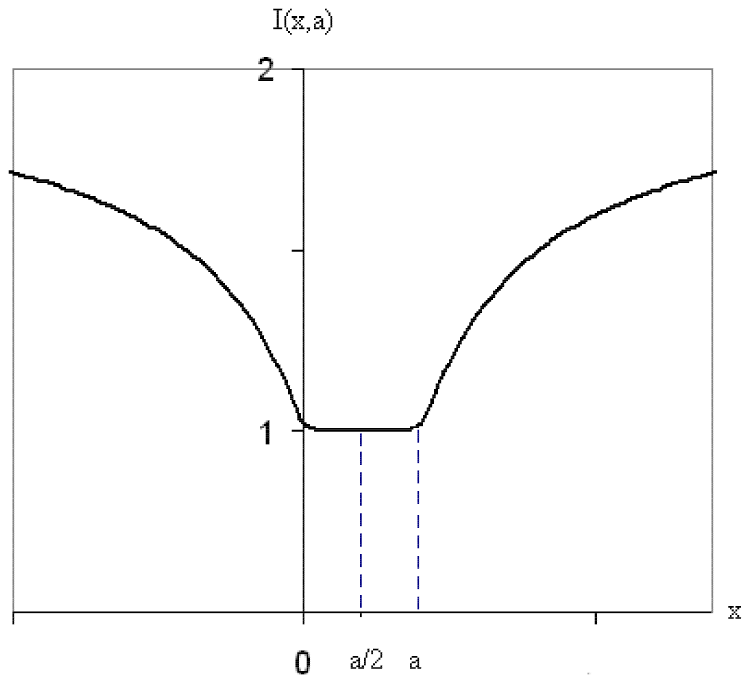}
\end{center}
\caption{Graph of ${\mathcal I}(\Gamma(x, a))$}
\label{inconstancygraph}
\end{figure}

\begin{figure}[!h]
\begin{center}
        \epsfig{file=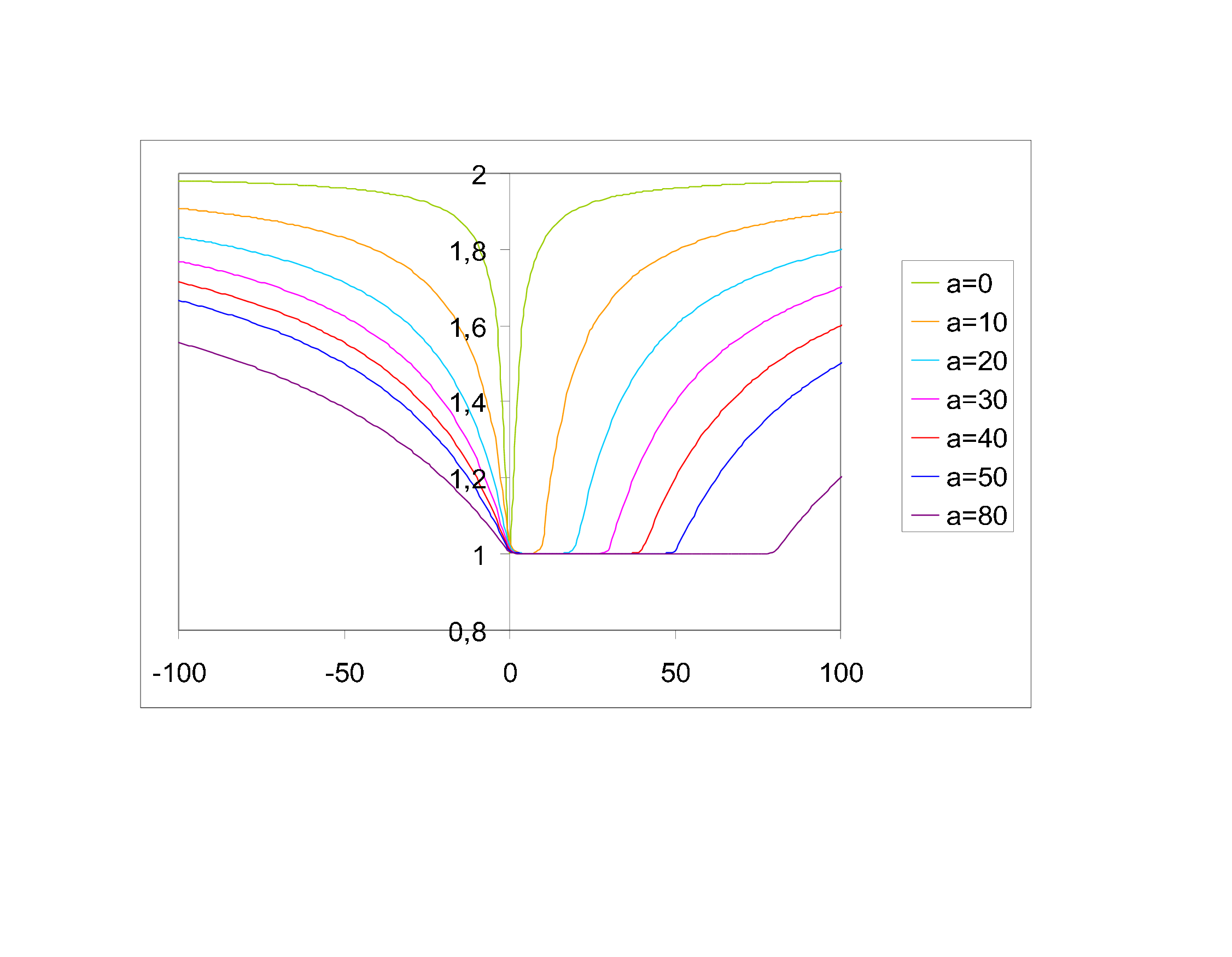}
\end{center}
\caption{Varying $a$ in the graph of ${\mathcal I}(\Gamma(x, a))$}
\label{varying-a}
\end{figure}

It is easy to show that ${\mathcal I}(\Gamma(a/2, a)) = 1$ (no fluctuation) and
$\lim_{x \to +\infty} {\mathcal I}(\Gamma(x, a)) = 2$ (when $x$ is large, the value
of $x$ is not really important, the inconstancy is close to $2$). We also have that
$({\mathcal I}(\Gamma(x, a)))' = 0$ if and only if $x = a/2$.
In particular the graph of the function ${\mathcal I}(\Gamma(x, a))$ has the
aspect shown in Figures~\ref{inconstancygraph} and \ref{varying-a}.

\medskip

\begin{figure}[h]
\begin{center}   
        \epsfig{file=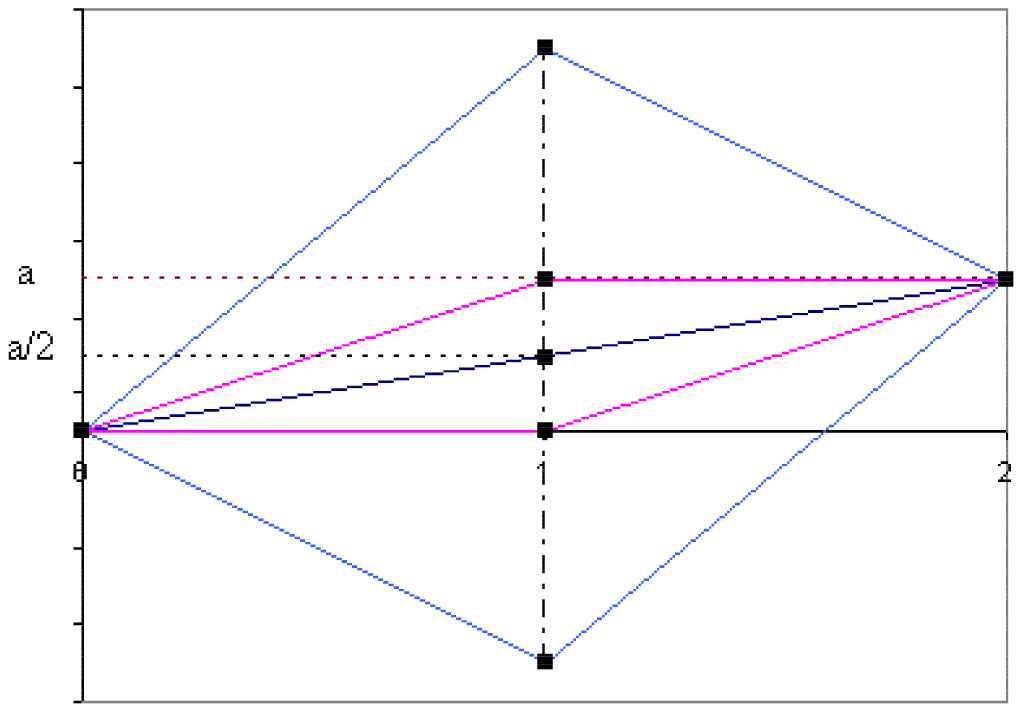}
\end{center}
\caption{Signs of slopes}
\label{slopes}
\end{figure}

We note that the curve is ``flat'' in the neighborhood of $a/2$,
or even for $x \in (0,a)$; see Figure~\ref{inconstancygraph}. 
This means that the inconstancy $({\mathcal I}(\Gamma(x, a)))$, 
which is equal to $1$ for $x=a/2$, remains close to $1$ when the two 
slopes of the curve have the same sign, while it is larger when the signs 
of the slopes are opposite, which correctly describes what a fluctuation 
should be (the MSE does not have this property); see Figure~\ref{slopes}. 
Also the inconstancy $({\mathcal I}(\Gamma(x, a)))$ tends quickly to $2$ 
when $a$ is small: see Figure~\ref{varying-a}. 

\subsection{General remarks}

If we look more generally at the inconstancy of
$\Gamma(a_1, a_2, \ldots, a_n)$, what will clearly matter for its size 
is the sequence of slopes: growth and signs of consecutive terms are
crucial characteristics of the sequence, which corresponds to the intuitive
idea of ``fluctuation''. Of course we always have the straightforward bounds
$$
1 \leq {\mathcal I}(\Gamma(a_1, a_2, \ldots, a_n)) \leq n
$$
(count the possible number of intersection points of 
$\Gamma(a_1, a_2, \ldots, a_n)$ with a random straight 
line and apply Theorem~\ref{Crofton}).

\bigskip

Conversely the inconstancy may be used to discriminate between curves,
i.e., to decide whether a curve fluctuates more than another, when
the ``visual aspect'' does not suffice to assert an intuitive answer.
We give two examples.

\subsection{Fluctuations of curves with four points}\label{nonintuitive}

\begin{figure}[!h]
\begin{center}
        \epsfig{file=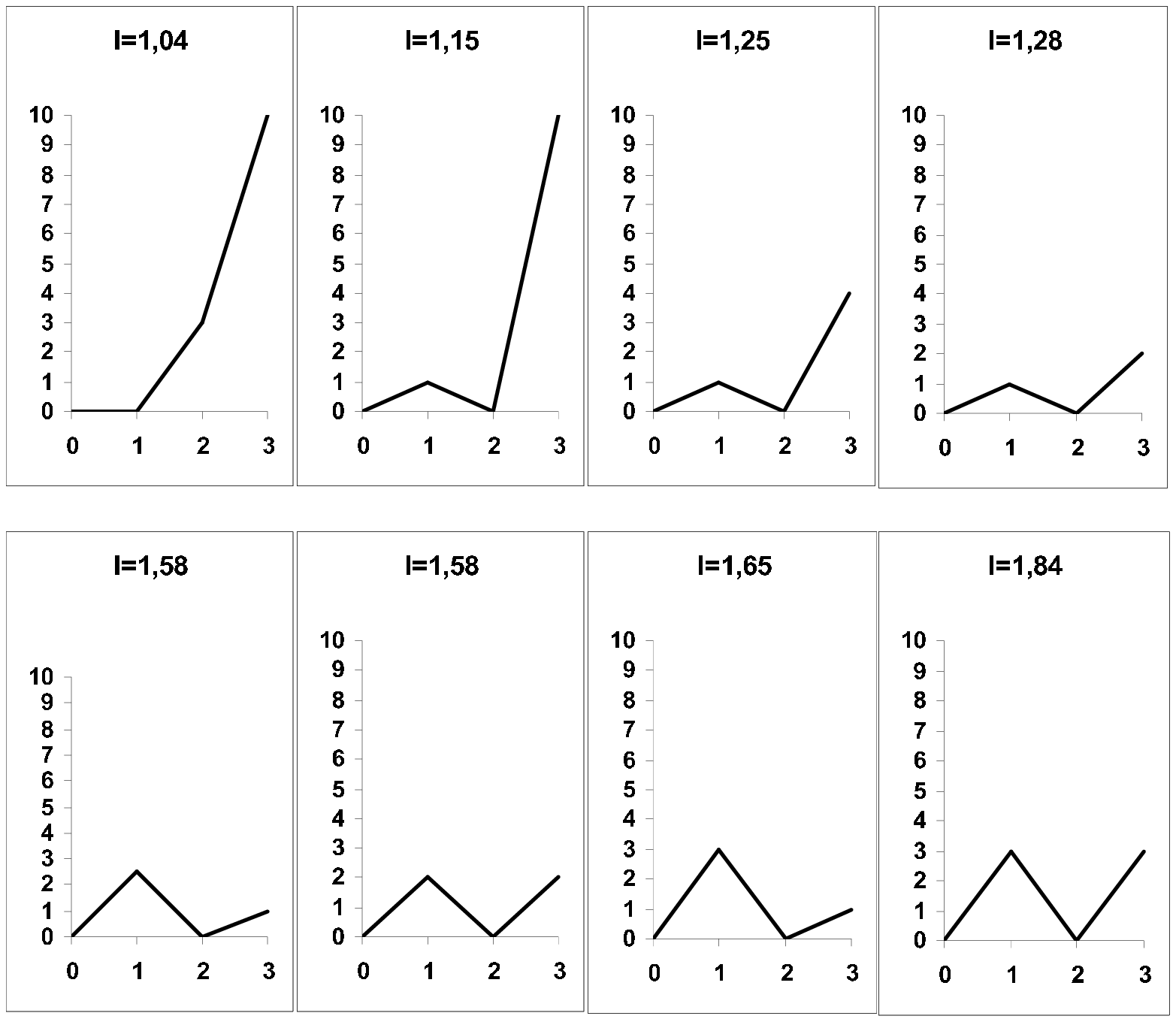}
\end{center}
\caption{Inconstancy discriminates between fluctuations}
\label{deciding}
\end{figure}

In Figure~\ref{deciding} inconstancies permit to discriminate between
``less fluctuating'' and ``more fluctuating'' curves, though there is
no visual evidence of which curve fluctuates more. It is interesting
to note that the maximum of the function is not really taken into account,
only the variations count (look, e.g.,  at the two examples with inconstancy 
$1.58$ in Figure~\ref{deciding}).

\subsection{A case where inconstancy does not discriminate}

\begin{figure}[!h]
\begin{center}
        \epsfig{file=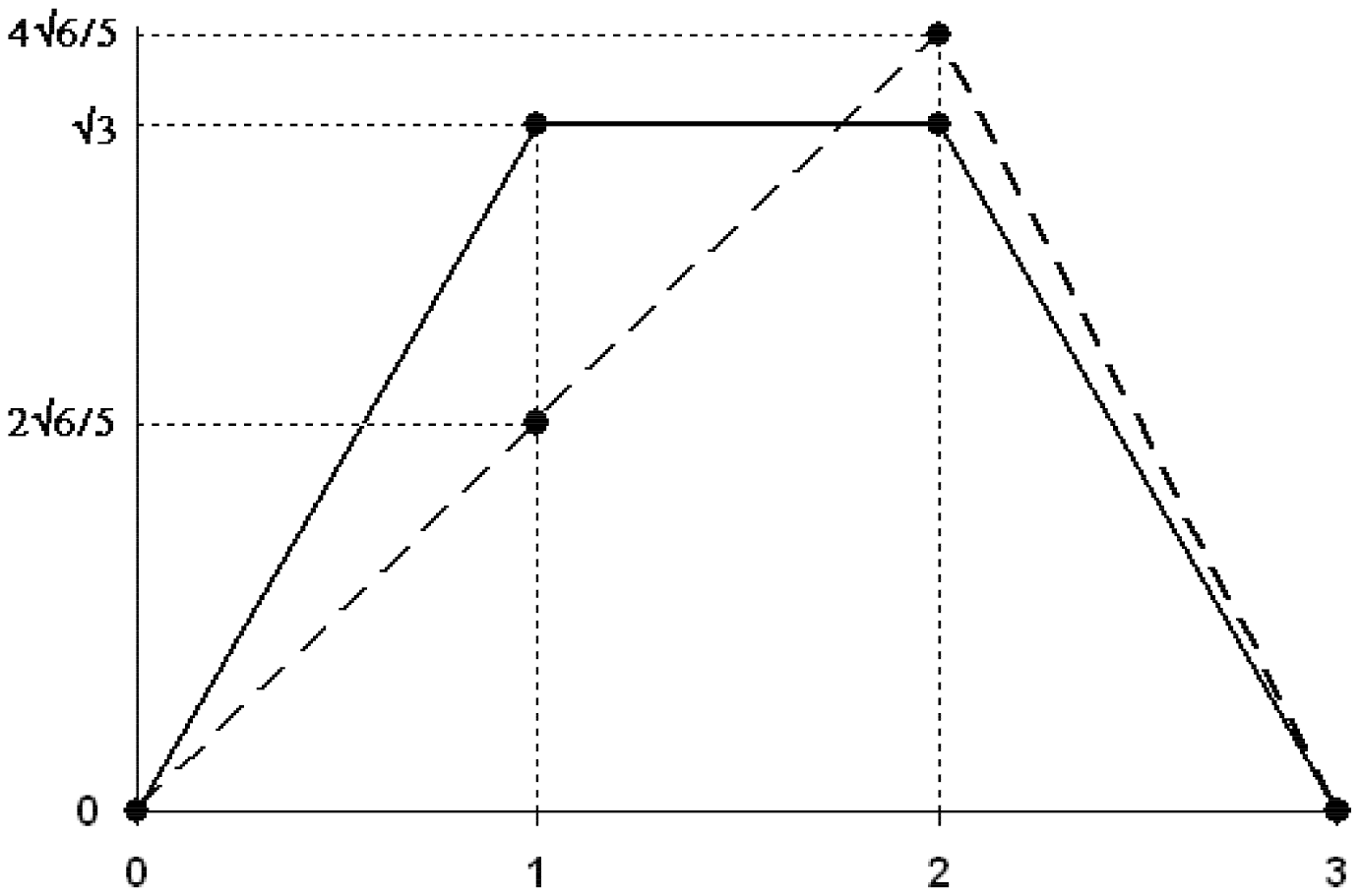}
\end{center}
\caption{Same inconstancy}
\label{sameinco}
\end{figure}

The lengths and inconstancies of the two curves $\Gamma(\sqrt{3}, \sqrt{3}, 0)$ 
and $\Gamma(2\sqrt{6}/5, 4\sqrt{6}/5, 0)$ (see Figure~\ref{sameinco}) are the same.

\section{Inconstancy of sequences}

Inconstancy of (finite or infinite) sequences can be defined in a 
straightforward way from what precedes. 

\begin{definition}
Let $(u_n)_{0 \leq n \leq N}$ be a finite sequence of real numbers, with 
$u_0 = 0$ say. Let $\Gamma_n$ be the union of the straight line segments 
$(0,0)$---$(1,u_1)$, $(1,u_1)$---$(2,u_2)$,...,$(n-1,u_{n-1})$---$(N,u_N)$,
then the inconstancy of $(u_n)_{0 \leq n \leq N}$ is defined by
$$
{\mathcal I}((u_n)_{0 \leq n \leq N}) := {\mathcal I}(\Gamma_N).
$$
Let $(u_n)_{n \geq 0}$ be an infinite sequence of real numbers, with
$u_0 = 0$ say. Then the inconstancy of $(u_n)_{n \geq 0}$ is defined by
$$
{\mathcal I}((u_n)_{n \geq 0}) := 
\limsup_{N \to \infty}{\mathcal I}((u_n)_{0 \leq n \leq N}) \ \ \ \
(\mbox{\rm or } \lim_{N \to \infty}{\mathcal I}((u_n)_{0 \leq n \leq N}) \
\mbox{\rm if the limit exists}). 
$$
\end{definition}

The inconstancy of an infinite sequence depends in particular of how long and 
frequently the sequence levels off: this is particularly clear for binary sequences 
as shown in Theorem~\ref{binary} below.

\begin{theorem}\label{binary}
\begin{itemize}

\item{(i)}
Let $(u_n)_{0 \leq n \leq N}$ be a finite sequence taking two values $0$ and 
$h > 0$, with $u_0 = 0$. Let $\alpha \geq 1$ be the index such that
$u_0 = u_1 = \ldots = u_{\alpha - 1} = 0$ and $u_{\alpha} \neq 0$.
In other words $\alpha$ is the length of the longest initial string of $0$'s.
Analogously let $\beta$ be the length of the longest final string of $0$'s.
If $\beta=0$, let $\gamma \geq 0$ be the largest index such that $u_{\gamma}=0$. 
Let ${\mathcal N}_{00}, {\mathcal N}_{hh}, {\mathcal N}_{0h}, {\mathcal N}_{h0}$ 
be respectively the number of blocks of the form $00$, $hh$, $0h$, $h0$ in the 
sequence. Then
$$
{\mathcal I}((u_n)_{0 \leq n \leq N}) =
\left\{
\begin{array}{ll}
\displaystyle
2 \frac{{\mathcal N}_{00} + {\mathcal N}_{hh}
+ (\sqrt{1 + h^2})({\mathcal N}_{0h} + {\mathcal N}_{h0})}
{\sqrt{h^2 + \alpha^2} + N - \alpha - \beta + \sqrt{h^2 + \beta^2} + N} \
&\mbox{\rm if } \beta > 0; \\
\displaystyle
2 \frac{{\mathcal N}_{00} + {\mathcal N}_{hh}
+ (\sqrt{1 + h^2})({\mathcal N}_{0h} + {\mathcal N}_{h0})}
{\sqrt{h^2 + \alpha^2} + N - \alpha  + \sqrt{h^2 + (N-\gamma)^2} + \gamma} \
&\mbox{\rm if } \beta = 0.
\end{array}
\right.
$$

\item{(ii)}
Let $(u_n)_{n \geq 0}$ be an infinite sequence taking two values $0$ and
$h > 0$, with $u_0 = 0$. We make the assumption that the frequencies
of occurrences of the blocks $00$, $hh$, $0h$, $h0$ in the sequence
exist and are respectively equal to ${\mathcal F}_{00}, {\mathcal F}_{hh}, 
{\mathcal F}_{0h}, {\mathcal F}_{h0}$. Then
$$
{\mathcal I}((u_n)_{n \geq 0}) =
{\mathcal F}_{00} + {\mathcal F}_{hh} + (\sqrt{1 + h^2})({\mathcal F}_{0h} 
+ {\mathcal F}_{h0}) = 
1 + (\sqrt{1 + h^2} - 1)({\mathcal F}_{0h} + {\mathcal F}_{h0}).
$$

Similarly let $(u_n)_{n \geq 0}$ be an infinite sequence taking only finitely 
many real values, and let $H$ be this set of values. We make the assumption that
the frequencies of occurrences of all length-$2$ blocks $jj'$ ($j, j' \in H$)
exist and are respectively equal to ${\mathcal F}_{jj'}$. Then
$$
\begin{array}{lll}
\displaystyle{\mathcal I}((u_n)_{n \geq 0}) &=&
\displaystyle\sum_{j \in H}{\mathcal F}_{jj} + 
\displaystyle\sum_{j,j' \in H,\ j < j'}(\sqrt{1 + (j'-j)^2})
({\mathcal F}_{jj'} + {\mathcal F}_{j'j}) \\
&=& 1 + \displaystyle\sum_{j,j' \in H,\ j < j'}
(\sqrt{1 + (j'-j)^2} - 1)({\mathcal F}_{jj'} + {\mathcal F}_{j'j}).
\end{array}
$$
\end{itemize}
\end{theorem}

\proof First let $(u_n)_{0 \leq n \leq N}$ be a finite sequence taking two 
values $0$ and $h > 0$. Let $\alpha \geq 1$ be the index such that
$u_0 = u_1 = \ldots = u_{\alpha - 1} = 0$ and $u_{\alpha} \neq 0$.
In other words $\alpha$ is the length of the longest initial string of $0$'s.
Analogously let $\beta \geq 0$ be the length of the longest final string of $0$'s.
Finally, if $\beta=0$, let $\gamma \geq 0$ be the largest index such that 
$u_{\gamma}=0$. It is almost immediate that the convex hull of the curve $\Gamma_N$ 
consists of the four straight line segments ($(0,0)$---$(\alpha,h)$), 
($(\alpha,h)$---$(N-\beta,h)$),
($(N-\beta,h)$---$(N,0)$), ($(0,0)$---$(N,0)$) if $\beta > 0$, and
($(0,0)$---$(\alpha,h)$), ($(\alpha,h)$---$(N,h)$), ($(0,0)$---$(\gamma,0)$),
($(\gamma,0)$---$(N,h)$) if $\beta = 0$ (there are only three segments if
$\gamma=0$, which implies $\alpha=1$). Hence 
$$
\delta(\Gamma_N) = 
\left\{
\begin{array}{ll}
\sqrt{h^2 + \alpha^2} + N - \alpha - \beta + \sqrt{h^2 + \beta^2} + N \
&\mbox{\rm if } \beta > 0; \\
\sqrt{h^2 + \alpha^2} + N - \alpha  + \sqrt{h^2 + (N-\gamma)^2} + \gamma \
&\mbox{\rm if } \beta = 0. 
\end{array}
\right.
$$
while the length of the curve is
$$
\ell(\Gamma_N) = {\mathcal N}_{00} + {\mathcal N}_{hh} 
+ (\sqrt{1 + h^2})({\mathcal N}_{0h} + {\mathcal N}_{h0}).
$$
This gives the first part of the theorem, namely
$$
{\mathcal I}((u_n)_{0 \leq n \leq N}) =
\left\{
\begin{array}{ll}
\displaystyle
2 \frac{{\mathcal N}_{00} + {\mathcal N}_{hh}
+ (\sqrt{1 + h^2})({\mathcal N}_{0h} + {\mathcal N}_{h0})}
{\sqrt{h^2 + \alpha^2} + N - \alpha - \beta + \sqrt{h^2 + \beta^2} + N} \
&\mbox{\rm if } \beta > 0; \\
\displaystyle
2 \frac{{\mathcal N}_{00} + {\mathcal N}_{hh}
+ (\sqrt{1 + h^2})({\mathcal N}_{0h} + {\mathcal N}_{h0})}
{\sqrt{h^2 + \alpha^2} + N - \alpha  + \sqrt{h^2 + (N-\gamma)^2} + \gamma} \
&\mbox{\rm if } \beta = 0.
\end{array}
\right.
$$

\medskip

In order to prove the second part of the theorem, we will directly address the 
case of a sequence $(u_n)_{n \geq 0}$ taking any finite number of values (the 
proof is simpler than our original one, thanks to a remark of one of the referees).
The length of the curve $\Gamma_n$ clearly is
$$
\sum_{j \in H} {\mathcal N}_{jj} + 
\sum_{j, j' \in H, \ j < j'} 
\sqrt{1+(j'-j)^2} ({\mathcal N}_{jj'} + {\mathcal N}_{j'j}).
$$
The perimeter of the convex hull of $\Gamma_n$ satisfies
$$
2 \sqrt{N^2+u_N^2} \leq \delta(\Gamma_n) 
\leq 2N + M_N, \ \mbox{\rm where} \ M_N := \max\{u_n, \ 0 \leq n \leq N\}
$$
(the inequality on the left is due to the fact that the perimeter is larger than twice
the distance between $(0,0)$ and $(N,u_N)$; the right inequality comes from the fact
that the length of the convex hull is less than the perimeter of the rectangle
$(0,0)$---$(0,M_N)$---$(N,M_N)$---$(N,0)$).
Since the sequence $(u_n)_{n \geq 0}$ takes only finitely many values, this shows that
$$
\delta(\gamma_N) = 2N + O(1).
$$
Hence
$$
{\mathcal I}((u_n)_{n \geq 0}) = \sum_{j \in H} {\mathcal F}_{jj} + 
\sum_{j, j' \in H, \ j < j'} 
\sqrt{1+(j'-j)^2} ({\mathcal F}_{jj'} + {\mathcal F}_{j'j}).
$$
But
$$
\sum_{j \in H} {\mathcal F}_{jj} +
\sum_{j, j' \in H, \ j < j'} ({\mathcal F}_{jj'} + {\mathcal F}_{j'j}) = 1
$$
hence the result.  \endpf

\section{Computing the inconstancy of classical sequences}

Using Theorem~\ref{binary} we see that the inconstancy of an infinite binary 
$(0,1)$-sequence must belong to the interval $[1, \sqrt{2}]$. Bound $1$ is reached 
if and only if
${\cal F}_{01} + {\cal F}_{10} = 0$, i.e., ${\cal F}_{01} = {\cal F}_{10} = 0$. 
Bound $\sqrt{2}$ is reached when ${\cal F}_{01} + {\cal F}_{10}= 1$, i.e.,
 ${\cal F}_{00} + {\cal F}_{11} = 0$, i.e., ${\cal F}_{00} = {\cal F}_{11} = 0$.
Let us compute the inconstancy of some classical binary sequences.

\subsection{Periodic sequences}

The sequence $(0^d1)^{\infty} = (00...01)^{\infty}$ (periodic of 
period $(d+1)$, where the period pattern consists of $d$ symbols $0$ followed 
by one symbol $1$). It is easy to compute ${\mathcal F}_{00} = \frac{d-1}{d+1}$, 
${\mathcal F}_{11} = 0$, ${\mathcal F}_{01} = {\mathcal F}_{10} = \frac{1}{d+1}$. Hence
$$
{\mathcal I}((0^d1)^{\infty}) = \frac{d-1 + 2\sqrt{2}}{d+1}.
$$
In particular, ${\mathcal I}((01)^{\infty}) = \sqrt{2} = 1.414...$
while ${\mathcal I}((0^d1)^{\infty})$ tends to $1$ when $d$ tends to infinity:
this corresponds to the fact that the curve becomes more and more flat when
$d$ increases. The case $d=1$ is somehow the worst case among periodic and
nonperiodic binary sequences in terms of levelling off (or flatness).

\subsection{Random sequences}

A random sequence of $0$'s and $1$'s. For almost all binary sequences 
we have ${\mathcal F}_{00} = {\mathcal F}_{11} = {\mathcal F}_{01} = {\mathcal F}_{10} 
= \frac{1}{4}$.
Hence if $(r_n)_{n \geq 0}$ is ``a random sequence'' of $0$'s and $1$'s, then
$$
{\mathcal I}((r_n)_{n \geq 0}) = \frac{1 + \sqrt{2}}{2} = 1.207...
$$

\subsection{Some automatic sequences}

We first recall a few notions of combinatorics on words; see, e.g. \cite{AS2}. 
A finite set is called an {\it alphabet}. Its elements are called {\it letters}.
For an alphabet $A$, we let $A^*$ denote the free monoid spanned by $A$ and 
equipped with the {\it concatenation}. Elements of $A^*$ are called {\it words} 
on $A$; the {\it length} of the word $a_1 a_2 ... a_n$, with $a_i \in A$, is 
$n$. Homomorphisms of monoids are called {\it morphisms}. A morphism from $A^*$ 
to $B^*$ is determined by the images of the letters in $A$. It is called
{\it uniform} if the images of all letters have the same length. The transition
matrix of a morphism $\sigma: A^* \to B^*$ counts the number of times
the letter $b_j$ in $B$ occurs in $\sigma(a_i)$. Finally a sequence is called
{\it automatic} if it is the pointwise image of a fixed point of a nontrivial
uniform morphism.

\bigskip

-- Recall that the {\em Thue-Morse sequence\,} with values $0$ and $1$
can be defined as the fixed point beginning with $0$ of the morphism $0 \to 01$,
$1 \to 10$ (see, e.g., \cite{AS1}): it is the most famous example of
{\em automatic sequences} (see, e.g., \cite{AS2}). The first few terms of the
Thue-Morse sequence $(m_n)_{n \geq 0}$ are
$$
0 \ 1 \ 1 \ 0 \ 1 \ 0 \ 0 \ 1 \ 1 \ 0 \ 0 \ 1 \ 0 \ 1 \ 1 \ 0 \ 1 \ 0 \ \ldots
$$
The frequencies of occurrences of blocks of length $2$ are given by
${\mathcal F}_{00} = {\mathcal F}_{11} = \frac{1}{6}$ and 
${\mathcal F}_{01} = {\mathcal F}_{10} = \frac{1}{3}$: 
this is a classical exercise that involves the morphism on four letters
defined by $a \to ab$, $b \to ca$, $c \to cd$, $d \to ac$. An alternative
proof consists of noting that the sequence $((m_n + m_{n+1}) \bmod 2)_{n \geq 0}$ 
is the {\em period doubling sequence}, i.e., the fixed point of the morphism
$1 \to 10$, $0 \to 11$; the sum of frequencies of the blocks $01$ and $10$
in the Thue-Morse sequence is thus the frequency of $1$'s in the period doubling
sequence which is easily seen to be $2/3$ (look at the transition matrix of
the morphism $1 \to 10$, $0 \to 11$). Hence
$$
{\mathcal I}((m_n)_{n \geq 0}) = \frac{1 + 2\sqrt{2}}{3} = 1.276...
$$
Note that the ``high'' value of this inconstancy is related to the absence
of long strings of $0$'s or of $1$'s: namely the Thue-Morse sequence does not 
contain the blocks $000$ and $111$. 

\bigskip

-- The {\em Shapiro-Rudin sequence\,} $(r_n)_{n \geq 0}$ with values $0$ and 
$1$ can be defined as the sequence of parities of the number of (possibly 
overlapping) $11$'s in he binary expansions of the integers 
$0, 1, 2, \ldots, n \ldots$ (see, e.g., \cite{AS2}). It is clear that the
sum of frequencies of occurrences of the blocks $01$ and $10$ is the frequency
of occurrences of the letter $1$ in the sequence $(r'_n)_{n \geq 0}$ defined
by $r'_n := (r_n + r_{n+1}) \bmod 2$. This last sequence is easily seen to be
the pointwise image under the map $a \to 0$, $b \to 0$, $c \to 1$, $d \to 1$
of the infinite fixed point of the morphism $a \to ab$, $b \to cd$, $c \to ad$,
$d \to cb$. (Hint: prove that both this pointwise image and the sequence
$(r'_n)_{n \geq 0}$ satisfy the recursive relations $r'_{4n} = 0$, 
$r'_{4n+1} = r'_{2n}$, $r'_{4n+2} = 1$, $r'_{4n+3} = 1 + r'_{2n+1} \bmod 2$, with 
$r'_0 = 0$. From this it is straightforward that the frequency of occurrences 
of $1$ in the sequence $(r'_n)_{n \geq 0}$ is equal to $1/2$. Hence
$$
{\mathcal I}((r_n)_{n \geq 0}) = \frac{1 + \sqrt{2}}{2} = 1.207...
$$
which is the same inconstancy as for a random sequence.

\bigskip

-- The {\em (regular) paperfolding sequence\,} $(z_n)_{n \geq 0}$ with values 
$0$ and $1$ can be defined by $z_{4n} = 0$, $z_{4n+1} = 1$, $z_{2n+1} = z_n$.
Reasoning as for the Shapiro-Rudin sequence (left to the reader) leads to
$$
{\mathcal I}((r_n)_{n \geq 0}) = \frac{1 + \sqrt{2}}{2} = 1.207...
$$
which is again the same inconstancy as for a random sequence.

\subsection{Sturmian sequences}

Recall that a Sturmian sequence can be defined as a (binary) sequence having
exactly $n+1$ blocks of length $n$ for every integer $n \geq 1$ (see, e.g.,
\cite{AS2, Lot}). In particular Sturmian sequences are not ultimately periodic,
and the blocks $00$ and $11$ cannot both occur in a same Sturmian sequence.
Since interchanging $0$'s and $1$'s in a Sturmian sequence gives a Sturmian 
sequence, we may suppose that no $11$ occurs. But then the frequencies of
occurrences of the blocks $01$ and $10$ in the sequence are both equal to
the frequency of occurrence of $1$, hence to the {\it slope\,} of the
Sturmian sequence (see, e.g., \cite[Theorem~10.5.8, page~318]{AS2}). Thus the 
inconstancy of a Sturmian sequence of slope $\alpha \in (0,1)$ without 
the block $11$ in it (resp.\ of slope $1 - \alpha \in (0,1)$ without
the block $00$ in it) is
$$
{\mathcal I} = 1 + 2(\sqrt{2} -1)\alpha.
$$
Recall that if the sequence does not contain the block $11$, then $\alpha$ 
belongs to $(0,1/2)$, hence as expected ${\mathcal I}$ belongs to $(1,\sqrt{2})$.

\begin{remark}
A possible application of inconstancy of infinite sequences can be to ``predict''
the $n$th term of a very long (or infinite) sequence knowing its first $n-1$ terms:
if $n$ is large enough, $u_n$ ``should'' be close to a value minimizing the difference 
$|{\mathcal I}(\Gamma_n) - {\mathcal I}(\Gamma_{n-1})|$.
\end{remark}

\begin{remark}\label{folding}
A different way of defining the inconstancy of a binary sequence could be to interpret
it as a sequence on the alphabet $\{$L(eft), R(ight)$\}$. Then to associate with this 
(LR) sequence a 2D curve drawn on the lattice ${\mathbb Z}^2$, consisting of horizontal 
and vertical segments. The first segment is $(0,0)$---$(1,0)$; then for each value of the LR
sequence we make a $\pm \pi/2$ turn. The inconstancy of the sequence could be defined
as the inconstancy of the curve obtained that way. The reader will have recognized curves
studied, e.g., in \cite{MFT}, where paperfolding sequences enter the picture. This notion
of inconstancy for sequences would thus be terminologically closer to the ``folding index''
of \cite{CYW}. Since the choice of $\pm \pi/2$ is arbitrary (another angle could have
been chosen), it is not clear whether this definition is pertinent or if one should
consider {\bf all} possible angles, thus obtaining a set of inconstancies for any
given sequence.
\end{remark}

\section{Algorithmic aspects}

In order to compute the inconstancy 
${\mathcal I}(\Gamma) := \frac{2 \ell(\Gamma)}{\delta(\Gamma)}$
of a piecewise affine curve $\Gamma$, the perimeter of the convex hull of 
$\Gamma$ is needed. Hence we need to construct the convex hull of a finite 
set consisting of, say, $n$ points. Several algorithms are available, 
their complexity is in ${\mathcal O}(n \log n)$ (see, e.g., the Graham scan 
studied in \cite{Graham}, the Jarvis march studied in \cite{Jarvis};
see also, e.g., the papers \cite{PH, Preparata} --in particular
\cite{Preparata} gives an optimal real-time algorithm for planar convex hulls).

\medskip

Implementations of these algorithms are classical in usual softwares:
for example the command {\tt convhull} in {\tt Maple} (with the package 
{\tt Convex}), the command {\tt ConvexHull} in {\tt Mathematica}, the 
command {\tt convex$\_$hull} in {\tt Scilab}, or the command {\tt convhull}
(see also {\tt convhulln}) in {\tt Matlab}. Also note that {\tt Qhull}
computes convex hulls, Delaunay triangulations, Voronoi diagrams, halfspace 
intersections about a point, furthest-site Delaunay triangulations, and 
furthest-site Voronoi diagrams (see {\tt http://www.qhull.org/}).
Demonstrations of computations can be found on several sites; see e.g.,

{\tt http://www.piler.com/convexhull/}
 
{\tt http://www.cs.princeton.edu/courses/archive/fall08/cos226/demo/ah/GrahamScan.html}

{http://www.cs.princeton.edu/courses/archive/fall08/cos226/demo/ah/ConvexHull.html}

{\tt http://www.cse.unsw.edu.au/$\sim$lambert/java/3d/hull.html}

\section{Conclusion}
Inspired by the theorem of Cauchy-Crofton, the inconstancy of a curve could be 
a way of detecting large fluctuations of a curve, different from (and hopefully
better than) usual indexes such as the residual variance. We intend to test
this idea in three domains: fluctuations of biological parameters \cite{VMAH},
fluctuations of the stockmarket \cite{AMstock} and smoothness of musical themes 
\cite{AMmus}. Two other directions could be the following. First, a way of 
discriminating between models that describe a given phenomenon with the same error 
bound (e.g., prediction of electric load and consumption) could be to choose the model 
for which the difference between data and predictions has maximal inconstancy (when the 
inconstancy is close to $1$, this difference is ``quasi-affine''; this means that there
is a ``quasi-affine'' bias in the model that can/should be corrected {\em a priori}).  
Second, we alluded to fractal-like ``chaotic'' (disordered) curves in the introduction; 
coming across, e.g., the paper \cite{Buttenfield} we recall that measuring the 
``complexity'' of geographic objects classically involves their fractal dimension and,
e.g., their ``length''; we could also think 
of looking at their inconstancy (typically how complicated a river can be, i.e., 
how far from straight it looks, can be measured by the number of intersection points  
with a random straight line). A natural question then occurs: to what extent fractal 
dimension and inconstancy are related? Or what can be said of the intersection with 
straight lines of a set with given fractal dimension? Such questions also make sense
for (in)finite sequences, in particular in view of Remark~\ref{folding}. Of course
the length of such curves is usually infinite while the length of the convex hull 
is finite (think of the von Koch curve). What could be looked at for fractals 
obtained by ``iteration'' is the inconstancy at each finite step of the iteration: 
it is conceivable that the fractal dimension shows up, though this is not the case 
for the von Koch curve.
Some ideas about these questions can be found, e.g., in 
\cite{Steinhaus, Marstrand, AkiSch, MFT, MF89, MF91}, in particular in relation with 
the {\it entropy} of a curve, as discussed in several papers of Mend\`es France. We 
will conclude this paper with that notion of {\em entropy\,} for a plane curve. Let 
$p_n$ be the probability that a straight line cuts the plane curve $\Gamma$ in exactly 
$n$ points, then the theorem of Cauchy-Crofton says that
$$
\sum_{n \geq 1} np_n = \frac{2 \ell(\Gamma)}{\delta(\Gamma)}\cdot
$$
It is natural to define the {\em entropy\,} of $\Gamma$ by 
$$
H(\Gamma) := \sum_{n \geq 1} p_n \log \frac{1}{p_n}\cdot
$$
Now how large can this expression be? Define the set of sequences ${\cal P}$ by 
$$
{\cal P} := \left\{(p_n)_{n \geq 1}; \ p_n \geq 0, \ \sum_{n \geq 1} p_n = 1, \ 
\sum_{n \geq 1} n p_n = \frac{2 \ell(\Gamma)}{\delta(\Gamma)}\right\},
\ \ \mbox{\rm and let} \ \ 
H_{\max}(\Gamma) := \max_{\cal P} H(\Gamma).
$$
It can be proven (see \cite{DKMF} for details, also see \cite{MF84}) that
$$
H_{\max}(\Gamma) = \log \frac{2 \ell(\Gamma)}{\delta(\Gamma)} + 
                   \frac{\beta}{e^{\beta}-1},
$$
where $\beta := \log \displaystyle\frac{2 \ell(\Gamma)}{2 \ell(\Gamma) - \delta(\Gamma)}$
(the quantity $\beta$ can be seen as the inverse of the {\em temperature\,} of the curve).

\medskip

A modified definition is thus proposed in \cite{MF83}, namely
$$
H(\Gamma) := \log \frac{2 \ell(\Gamma)}{\delta(\Gamma)}\cdot
$$ 
This definition was used in several papers (see, e.g., \cite{GW, DenCre, BCC}). 
With our terminology, it reads, as noted by Mend\`es France, ``the entropy is the 
logarithm of the inconstancy''. 
The reader might think of comparing this statement with the classical Weber-Fechner 
law in psychophysics according to which ``sensation is proportional to the logarithm 
of excitation'' (\cite{Fechner}; see also
{\tt http://psychclassics.yorku.ca/Fechner/}).

\section{Acknowledgments} The authors would like to thank several readers of this 
preprint, in particular the referees and P. Duchet. They also wish to thank very 
warmly M. Mend\`es France for interesting discussions and for several suggestions 
after he read a previous version of the paper.


\begin{thebibliography}{99}

\bibitem{AMmus} C. Agon, J.-P. Allouche, M. Andreatta, L. Maillard-Teyssier, 
Smoothness of musical pieces and inconstancy of discrete curves, 
{\it in preparation}.

\bibitem{AkiSch} S. Akiyama, K. Scheicher, Intersecting two-dimensional fractals
with lines, {\it Acta Sci. Math.} {\bf 71} (2005) 555--580.

\bibitem{AMstock} J.-P. Allouche, L. Maillard-Teyssier, Another way of quantifying
the fluctuations of the stockmarket, {\it in preparation}.

\bibitem{AS1} J.-P.~Allouche, J.~Shallit, The ubiquitous Prouhet-Thue-Morse
sequence, in C.~Ding, T.~Helleseth, H.~Niederreiter (Eds.), {\em Sequences
and their applications, Proceedings of SETA'98}, Springer, 1999, pp.~1--16.

\bibitem{AS2} J.-P. Allouche, J. Shallit, {\em Automatic sequences. Theory,
applications, generalizations}, Cambridge University Press, Cambridge, 2003.

\bibitem{AD} S. Ayari, S. Dubuc, La formule de Cauchy sur la longueur d'une courbe,
{\it Canad. Math. Bull.} {\bf 40} (1997) 3--9.

\bibitem{BCC} A. Balestrino, A. Caiti, E. Crisostomi, Generalised entropy of curves 
for the analysis and classification of dynamical systems, {\it Entropy\,} {\bf 11}
(2009) 249--270.

\bibitem{BKT} W. Z. Billewicz, W. F. F. Kemsley, A. M. Thomson, Indices of adiposity,
{\it Br. J. Prev. Soc. Med.} {\bf 16} (1962) 183--188. 

\bibitem{Blaschke1} W. Blaschke, {\it Vorlesungen \"uber Integralgeometrie I},
Hamburger Math. Einzelschr. {\bf 20}, Leipzig, B. G. Teubner, 1935.

\bibitem{Blaschke2} W. Blaschke, {\it Vorlesungen \"uber Integralgeometrie II},
Hamburger Math. Einzelschr. {\bf 22}, Leipzig, B. G. Teubner, pp.\ 61-127, 1937.

\bibitem{Buffon} G. L. Leclerc Buffon, {\it Histoire naturelle g\'en\'erale et
particuli\`ere}, Suppl\'ement, Tome quatri\`eme, Imprimerie Royale, 1777.


\bibitem{Buttenfield} B. Buttenfield, Treatment of the cartographic line, 
{\it Cartographica\,} {\bf 22} (1985) 1--26.

\bibitem{Cauchy1} A. Cauchy, Notes sur divers th\'eor\`emes relatifs \`a la
rectification des courbes, et \`a la quadrature des surfaces, {\it C.\ R. \ Acad. 
\ Sci. \ Paris\,} {\bf 13} (1841) 1060--1063. (Also in {\it {\OE}uvres compl\`etes} 
{\bf 6}, Gauthier-Villars, Paris, pp. 369--375, 1888.)

\bibitem{Cauchy2} A. Cauchy, M\'emoire sur la rectification des courbes et la
quadrature des surfaces courbes, {\it M\'em. Acad. Sci. Paris\,} {\bf 22}
(1850) 3--15. (Also in {\it {\OE}uvres compl\`etes} {\bf 2}, Gauthier-Villars, 
Paris, pp. 167--177, 1908.)

\bibitem{CYW} Y. K. Choi, D. M. Yan, W. Wang, A folding index of 2D curves,
{\it Comput.-Aided Des. Appl.} {\bf 1} (2004) 741--749.

\bibitem{Cor1} P. Cordier, M. Mend\`es France, P. Bolon, J. Pailhous, Entropy, degrees 
of freedom and free climbing: a thermodynamic study of a complex behavior based on 
trajectory analysis, {\it Int. J. Sport Psychol.} {\bf 24} (1993) 370--378.

\bibitem{Cor2} P. Cordier, M. Mend\`es France, P. Bolon, J. Pailhous, Thermodynamics 
study of motion behavior optimization, {\it Acta Biotheoretica\,} {\bf 42} (1994) 187--201.

\bibitem{Cor3} P. Cordier, M. Mend\`es France, J. Pailhous, P. Bolon, Entropy as a 
global variable of learning process, {\it Human Mov. Sci.} {\bf 13} (1994) 745--763.

\bibitem{CMSPNZCA} M. Cremona, M. H. P. Mauricio, L. C. Scavarda Do Carmo, R. Prioli, 
V. B. Nunes, S. I. Zanette, A. O. Caride, M. P. Albuquerque, Grain size distribution 
analysis in polycrystalline LiF thin films by mathematical morphology techniques on 
AFM images and X-ray diffraction data, {\it J. Microscopy\,} {\bf 197} (2000) 260--267.

\bibitem{Crofton} M. W. Crofton, On the theory of local probability, applied to 
straight lines drawn at random in a plane; the methods used being also extended 
to the proof of certain new theorems in the Integral Calculus, {\it  Philos. 
Trans. R. Soc. Lond.} {\bf 158} (1868) 181--199.

\bibitem{DenCre} A. Denis, F. Cr\'emoux, Using the entropy of curves to segment
a time or spatial series, {\it Math. Geol.} {\bf 34} (2002) 899--914.

\bibitem{DKMF} Y. Dupain, T. Kamae, M. Mend\`es France, Can one measure the 
temperature of a curve?, {\it Arch. Rational Mech. Anal.} {\bf 94} (1986) 
155--163.  Corrigenda, {\it Arch. Rational Mech. Anal.} {\bf 98} (1987) 395.

\bibitem{Fechner} G. T. Fechner, {\it Elemente der Psychophysik}, Breitskopf und
H\"artel, Leipzig, 1860.

\bibitem{Furuyama} M. Furuyama, Applications of L-mosaic map to spatial 
analysis of urban map, Mem. Fac. Ind. Arts, Kyoto Technical University,
Science and Technology, {\bf 25} (1976) 113--134.

\bibitem{GS} N. Y. Gnedin, S. F. Shandarin, Morphology of the secondary cosmic 
microwave background anisotropies: the key to `smouldering' reionization,
{\it Monthly Notices Royal Astron. Soc.} {\bf 337} (2002) 1435--1440.

\bibitem{GW} G. Gouesbet, M. E. Weill, Complexities and entropies of periodic 
series with application to the transition to turbulence in the logistic map,
{\it Phys. Rev. A\,} {\bf 30} (1984) 1442--1448.

\bibitem{Graham} R. L. Graham, An efficient algorithm for determining the 
convex hull of a finite planar set, {\it Inf. Process. Lett.} {\bf 1} 
(1972) 132--133.

\bibitem{HCLB} Y. J. Han, Y.-J. Cho, W. E. Lambert, C. K. Bragg, Identification 
and measurement of convolutions in cotton fiber using image analysis,
{\it Art. Intel. Review\,} {\bf 12} (1998) 201--211.

\bibitem{Jarvis} R. A. Jarvis, On the identification of the convex hull of a 
finite set of points in the plane, {\it Inf. Proc. Letters} {\bf 2} (1973) 18--21. 

\bibitem{Kaiser} L. Kaiser, Unbiased estimation in line-intercept sampling, 
{\it Biometrics\,} {\bf 39} (1983) 965--976.

\bibitem{Laplace} P. S. de Laplace, {\it Th\'eorie Analytique des Probabilit\'es}, 
V. Courcier, 1812.

\bibitem{Langevin1} R. Langevin, {\it Introduction to integral geometry}, Col\'oquio
Brasileiro de Mathem\'atica. [21st Brazilian Mathematics Colloquium], Instituto de
Mathem\'atica Pura e Aplicada (IMPA), Rio de Janeiro, 1997.

\bibitem{Langevin2} R. Langevin, {\it Integral geometry from Buffon to the use of 
twentieth century mathematics}, 2009, preprint (or updated version) available at 
\newline
\noindent
{\tt http://math.u-bourgogne.fr/topolog/langevin/preprints.html}

\bibitem{Lot} M.~Lothaire, {\it Algebraic Combinatorics On
Words, Encyclopedia of Mathematics and its Applications}, vol.~90,
Cambridge University Press, 2002.

\bibitem{Marstrand} J. M. Marstrand, Some fundamental geometrical properties of plane 
sets of fractional dimensions, {\it Proc. Lond. Math. Soc., III. Ser.} {\bf 4} (1954)
257--302.

\bibitem{MF82} M. Mend\`es France, Paper folding, space-filling curves and Rudin-Shapiro 
sequences, {\it Contemp. Math.} {\it 9} (1982) 85--95.

\bibitem{MF83a} M. Mend\`es France, Chaotic curves, in {\it Rhythms in biology and
other fields of application}, Proc. Journ. Soc. Math. France, Luminy, 1981, 
Lect. Notes Biomath. {\bf 49} (1983) 352--367.

\bibitem{MF83} M. Mend\`es France, Les courbes chaotiques, {\it Images de la Physique (CNRS)},
{\bf 51} (1983) 5--9, available electronically at
\newline
{\tt http://www.cnrs.fr/publications/imagesdelaphysique/Archives-1975-1989/1983/1983-5-9.pdf}

\bibitem{MF84} M. Mend\`es France, Folding paper and thermodynamics, {\it Phys. Rep.}
{\bf 103} (1984) 161--172.

\bibitem{MF89} M. Mend\`es France, Chaos implies confusion, in {\it Number theory and 
dynamical systems}, Lond. Math. Soc. Lect. Note Ser. {\bf 134} (1989) 137--152. 

\bibitem{MF91} M. Mend\`es France, The Planck constant of a curve, in {\it Fractal geometry 
and analysis (Montreal, PQ, 1989)}, NATO Adv. Sci. Inst. Ser. C Math. Phys. Sci. {\bf 346},
Kluwer Acad. Publ., Dordrecht, 1991, pp.~325--366. 

\bibitem{MF06} M. Mend\`es France, Poincar\'e et les probabilit\'es g\'eom\'etriques,
in {\it L'h\'eritage scientifique d'Henri Poincar\'e}, sous la direction d'\'E. Charpentier,
\'E. Ghys, A. Lesne, Belin, 2006, ch.~15, pp.~316--330.

\bibitem{MFH} M. Mend\`es France, A. H\'enaut, Art, therefore entropy, {\it Leonardo\,} 
{\bf 27}, Art and Science Similarities, Differences and Interactions: Special Issue (1994) 
219--221.
 
\bibitem{MFT} M. Mend\`es France, G. Tenenbaum, Dimension des courbes planes, papiers 
pli\'es et suites de Rudin-Shapiro, {\it Bull. Soc. Math. France\,} {\bf 109} (1981)
207--215.

\bibitem{Myers} S. B. Myers, Review: W. Blaschke, Integralgeometrie, and 
L. A. Santalo, Integralgeometrie, and W. Blaschke, Vorlesungen \"uber 
Integralgeometrie, Vol. 1, and W. Blaschke, Vorlesungen über Integralgeometrie,
Vol. 2, and W. Blaschke, \"uber eine geometrische Frage von Euclid bis Heute,
{\it  Bull. Amer. Math. Soc.} {\bf 44} (1938) 614--615.

\bibitem{Nesetril1} J. Ne\v{s}et\v{r}il, Art of graph drawing and art,
{\it J. Graph Algorithms Appl.} {\bf 6} (2002) 131--147.

\bibitem{Nesetril2} J. Ne\v{s}et\v{r}il, Aesthetics for computers or how to measure
harmony, in {\it The visual mind, II}, M. Emmer ed., Leonardo Books, 2005, pp.~35--58.

\bibitem{Nicolis} J. S. Nicolis, Chaotic dynamics applied to information processing,
{\it Rep. Prog. Phys.} {\bf 49} (1986) 1109--1196.

\bibitem{PPVA} M. N. Pons, V. Plagnieux, H. Vivier, D. Audet, Comparison of methods 
for the characterisation by image analysis of crystalline agglomerates: the case of 
gibbsite, {\it Powder Technol.} {\bf 157} (2005) 57--66.  

\bibitem{Preparata} F. P. Preparata, An optimal real-time algorithm for planar 
convex hulls, {\it Comm. ACM} {\bf 22} (1979) 402--405.

\bibitem{PH} F. P. Preparata, S. J. Hong, Convex hulls of finite sets of points 
in two and three dimensions, {\it Comm. ACM} {\bf 20} (1977) 87--93.

\bibitem{Quetelet} A. Quetelet, {\it Sur l'homme et le d\'eveloppement de ses
facult\'es ou Essai de physique sociale}, tome second, Bachelier, Paris, 1835.

\bibitem{RC2005} M.-F. Rolland-Cachera, Rate of growth in early life: a 
predictor of later health?, in {\it Early Nutrition and its Later Consequences: 
New Opportunities, Perinatal Programming of Adult Health -- EC Supported Research},
Advances in Experimental Medicine and Biology {\bf 569}, N. Back, I. R. Cohen, 
D. Kritchevsky, A. Lajtha, R. Paoletti, B. Koletzko, P. Dodds, H. Akerblom, M. Ashwell,
eds., Springer, 2005, pp. 35--39.

\bibitem{RC84} M.-F. Rolland-Cachera, M. Deheeger, F. Bellisle, M. Semp\'e, 
M. Guilloud-Bataille, \'E. Patois, Adiposity rebound in children: a simple indicator 
for predicting obesity, {\it Am. J. Clin. Nutr.} {\bf 39} (1984) 129--135. 

\bibitem{RC82} M.-F. Rolland-Cachera, M. Semp\'e, M. Guilloud-Bataille, 
\'E. Patois, F. Pequignot-Guggenbuhl, V. Fautrad, Adiposity indices in children, 
{\it Am. J. Clin. Nutr.} {\bf 36} (1982) 178--184.

\bibitem{Santalo} L. A. Santal\'o, {\it Introduction to integral geometry},
Actualit\'es Sci. Ind., no. 1198, Publ. Inst. Math. Univ. Nancago II. Herman et Cie, 
Paris, 1953. 

\bibitem{SBMSSS} J. Schmalzing, T. Buchert, A. L. Melott, V. Sahni, 
B. S. Sathyaprakash, S. F. Shandarin, Disentangling the cosmic web I: 
morphology of isodensity contours, {\it Astrophys. J.} {\bf 526} (1999) 568--578.
 
\bibitem{Steinhaus} H. Steinhaus, Length, shape, and area, {\it Colloquium Math.} 
{\bf 3} (1954) 1--13.

\bibitem{Teissier} B. Teissier, Volumes des corps convexes, g\'eom\'etrie et alg\`ebre, 
Le\c{c}on 7 in {\it Le\c{c}ons de math\'ematiques d'aujourd'hui} {\bf 3}, sous la direction 
d'\'E. Charpentier et N. Nikolski, Cassini, 2007.

\bibitem{Vergnaud} A.-C. Vergnaud, S. Bertrais, J.-M. Oppert, L. Maillard-Teyssier, 
P. Galan, S. Hercberg, S. Czernichow, Weight fluctuations and risk for metabolic syndrome 
in an adult cohort, {\it Int. J. Obesity\,} {\bf 32} (2008) 315--321.

\bibitem{VMAH} A.-C. Vergnaud, J. Oudinet, L. Maillard-Teyssier, J.-P. Allouche, 
S. Hercberg, Inconstancy of discrete curves and fluctuations of biological parameters,
{\it in preparation}.

\bibitem{WVDGHKLMN} T. Wagner, A. Vetter, N. Dimovic, S. E. Guber, W. Helmberg,
W. Kr\"ull, G. Lanzer, W. R. Mayr, J. Neum\"uller, Ultrastructural changes and 
activation differences in platelet concentrates stored in plasma and additive 
solution, {\it Transfusion\,} {\bf 42} (2002) 719--727.

\bibitem{ZPN} K. Zhao, S. Popescu, R. Nelson, Quantifying variances of 
Line-Intercept-Sampling estimators of percentage cover, available at

\noindent
{\tt http://sslsnap02.tamu.edu/EXCHANGE/KZhao/Kaiguang$\_$dissertation/japanese$\_$manuscript/Zhao$\_$text.doc}

\noindent
[see also J. Forest Planning (Japanese Society of Forest Planning) {\bf 13} (2008) 195--205].


\end{thebibliography}
\end{document}